\documentclass{amsart}
\usepackage{amssymb}
\usepackage{amsmath}

\usepackage{a4}
\usepackage{color}

\newcommand{\pu}[1]{\partial_{u_{#1}}}
\newcommand{\pv}[1]{\partial_{v_{#1}}}
\newcommand{\ps}[1]{\partial_{s_{#1}}}
\newcommand{\ep}{\varepsilon}
\newcommand{\BB}{\mathcal{B}}
\newcommand{\CDspan}{\text{Span}}

\def\mf{\mathcal{M}_F}
\def\ms{\mathcal{V}}
\def\gv{\mathcal{G}_\ms}

\begin{document}

\newtheorem{theorem}{Theorem}[section]
\newtheorem{lemma}[theorem]{Lemma}
\newtheorem{remark}[theorem]{Remark}
\newtheorem{definition}[theorem]{Definition}
\newtheorem{corollary}[theorem]{Corollary}
\newtheorem{example}[theorem]{Example}
\newtheorem{problem}{Problem}[section]
\def\qedbox{\hbox{$\rlap{$\sqcap$}\sqcup$}}
\makeatletter
  \renewcommand{\theequation}{%
   \thesection.\alph{equation}}
  \@addtoreset{equation}{section}
 \makeatother
\title[Curvature Homogeneous Manifolds]
{A New Family of Curvature Homogeneous Pseudo-Riemannian Manifolds}
\author{Corey Dunn}

\begin{address}{CD: Mathematics Department, California State University at San Bernardino,
San Bernardino, CA 92407, USA. Email: \it
cmdunn@csusb.edu.}\end{address}
\begin{abstract} We construct a new family of curvature homogeneous pseudo-Riemannian
manifolds modeled on $\mathbb{R}^{3k+2}$ for integers $k \geq 1$. In contrast to
previously known examples, the signature may be chosen to be $(k+1+a, k+1+b)$ where
$a,b \in \mathbb{N} \bigcup \{0\}$ and $a+b = k$.  The structure group of the $0$-model
of this family is studied, and is shown to be indecomposable.  Several invariants that
are not of Weyl type are found which will show that, in general, the members of this
family are not locally homogeneous.\end{abstract}
 \maketitle

\section{Introduction} 

Let   $(M,g)$ be a smooth pseudo-Riemannian manifold of signature $(p,q)$, and let $P \in M$.   Using the Levi-Civita connection $\nabla$, one can compute the Riemann curvature tensor $R\in \otimes ^4 T^*_PM$ as follows:
$$
R(X,Y,Z,W) := g(\nabla_X \nabla_Y Z - \nabla_Y \nabla_X Z - \nabla_{[X,Y]} Z, W), \text{ for } X,Y,Z,W \in T_PM\,.
$$
One similarly defines the tensors $\nabla^i R$, for $ i = 0, 1, 2, \ldots$.  For convenience,   we write $\nabla^0 R = R$.
 Let $g_P$, $R_P$, and $\nabla^iR_P$ denote
the evaluation of these tensors at the point $P$.

The manifold $(M,g)$ is \emph{r-curvature homogeneous} if for all points $P,Q \in M$ and $i = 0, 1, \ldots, r$, 
there exists a linear isomorphism $\Phi_{PQ}:T_PM \to T_QM$ so that $\Phi_{PQ}^* g_Q = g_P$ and $\Phi_{PQ}^*\nabla^iR_Q = \nabla^i R_P$.   

There is an equivalent characterization of $r$-curvature homogeneous manifolds that will be of use.  Let $V$ be a finite dimensional real vector space, let the dual vector space $V^* := \text{Hom}_{\mathbb{R}}(V, \mathbb{R})$, and let $( \cdot, \cdot)$
 be a symmetric nondegenerate inner product on $V$.  An element $A^0 \in \otimes^4 V^*$ is called an \emph{algebraic curvature tensor} on $V$ if it satisfies the following three properties for all $v_1, \ldots, v_4 \in V$:
$$ \begin{array}{r c l}
A^0(v_1, v_2, v_3, v_4)& = & - A^0(v_2, v_1, v_3, v_4), \\
A^0(v_1, v_2, v_3, v_4) & = & A^0(v_3, v_4, v_1, v_2), \text{ and } \\
0 & = & A^0(v_1, v_2, v_3, v_4) + A^0(v_2, v_3, v_1, v_4) \\
 & &\quad\qquad+ A^0(v_3, v_1, v_2, v_4)\,.
\end{array}$$
 An element $A^1 \in \otimes^5 V^*$ is called an \emph{algebraic covariant derivative curvature tensor} on $V$ if it satisfies the following four properties for all $v_1, \ldots, v_5 \in V$:
 $$ \begin{array}{r c l}
A^1(v_1, v_2, v_3, v_4;v_5)& = & - A^1(v_2, v_1, v_3, v_4;v_5), \\
A^1(v_1, v_2, v_3, v_4;v_5) & = & A^1(v_3, v_4, v_1, v_2;v_5), \\
0 & = & A^1(v_1, v_2, v_3, v_4;v_5) + A^1(v_2, v_3, v_1, v_4;v_5) \\
 & &\quad\qquad+ A^1(v_3, v_1, v_2, v_4;v_5), \\
 0 & = & A^1(v_1, v_2, v_3, v_4;v_5) + A^1(v_1, v_2, v_4, v_5;v_1) \\
 & &\quad\qquad+ A^1(v_1, v_2, v_5, v_1;v_4)\,.
\end{array}$$
Let $A^i \in \otimes^{4+i}V^*$ for $i = 2, 3, \ldots, r.$ The tensors $A^0$ and $A^1$ are algebraic analogues of $R$ and $\nabla R$.  The symmetries of the tensors $\nabla^2 R$, $\nabla^3 R$, \ldots  \ are more difficult to express and are not relevant to our discussion.  Thus, we will not impose any restrictions on the tensors $A^i$ for $i = 2, 3, \ldots, r.$  We define an
\emph{$r$-model} to be a tuple
$\ms_r:=(V,
(\cdot,\cdot), A^0,
\ldots, A^r)$.  A
\emph{weak
$r$-model} is an $r$-model without the   bilinear form. Thus, a pseudo-Riemannian manifold $(M,g)$ is $r$-curvature
homogeneous if and only if for each $P\in M$ there exists a linear isometry $\Phi_P:T_PM \to V$, with $\Phi^*_P A^i = \nabla^i R_P$ for $i =
0, 1, \ldots, r$.  In such an event we say that $(M,g)$ is $r$-modeled on $\ms_r$, or that $\ms_r$ is a $r$-model for $(M,g)$.  The
\emph{structure group} ${\gv}_{,r}$ of the  $r$-model $\ms_r$ is the group of isomorphisms of $\ms_r$.  For an $r$-curvature homogeneous space, this group is independent of
$P$.  

    It is clear that a locally homogeneous manifold is $r$-curvature homogeneous for all $r$.  The converse, however, is not always true:  There exist pseudo-Riemannian manifolds which are $r$-curvature homogeneous for some $r$, and not (locally) homogeneous.   The study of curvature homogeneity in the Riemannian setting began with a paper by I.M. Singer \cite{S} in 1960.  His result was extended by Podesta and Spiro  \cite{PS04} to the pseudo-Riemannian setting in 1996: 
\begin{theorem} \label{theorem:singer}
Let $(M,g)$ be a smooth, simply connected, complete manifold of dimension $n$.
 \begin{enumerate}
 \item  (Singer, 1960) If $(M,g)$ is Riemannian, then there exists an integer $k_{0,n}$ so that if $(M,g)$ is $k_{0,n}$-curvature homogeneous, then it is homogeneous.
  \item  (Podesta, Spiro, 1996) If $(M,g)$ is a pseudo-Riemannian manifold of signature $(p,q)$, then there exists an integer $k_{p,q}$ so that if $(M,g)$ is $k_{p,q}$-curvature homogeneous, then it is homogeneous.
 \end{enumerate}
 \end{theorem}

Since then, many authors have studied curvature homogeneous  manifolds both in the Riemannian and higher signature settings--indeed, the list of references is becoming quite large and we only summarize the results pertinent to our goal--for more details see  \cite{BKV, GBk2}.  Opozda \cite{Op} has obtained a result similar to Theorem \ref{theorem:singer} in the affine case.

In the Riemannian setting, it is clear that $k_{0,2} = 0$, and the efforts of of Gromov \cite{Gr} and Yamato \cite{Y} have established bounds on $k_{0,n}$ which are linear in $n$.  The work of Sekigawa, Suga, and Vanhecke \cite{SSV1, SSV2} shows $k_{0, 3} = k_{0,4} = 1$.   There are examples of $0$-curvature homogeneous Riemannian manifolds which are not locally homogeneous, see \cite{FKM, KTV,  T}.  There are no known examples of $1$-curvature homogeneous Riemannian manifolds which are not locally homogeneous.

In the pseudo-Riemannian setting, the situation is somewhat similar.  There are many known examples of $0$-curvature homogeneous pseudo-Riemannian manifolds  which are not locally homogeneous, see for example  
\cite{Bue, GSL} in the Lorentzian setting, and \cite{DG, GS, GS04, G-S} in the higher signature setting.  It is clear that $k_{1,1} = 0$.  The work of
Bueken, and  Djori\'c  \cite{BD} and  the work of Bueken and Vanhecke \cite{BV} shows that $k_{1,2} \geq 2$, while the work in \cite{DGS} shows
$k_{2,2} \geq 2$.  Derdzinski \cite{Derd} has also studied isometry invariants in signature $(2,2)$.  In contrast to the Riemannian setting, however, 
there exist examples of higher curvature homogeneity in the higher signature setting.  For instance,  examples constructed by  Gilkey and Nik\v
cevi\'c  \cite{GS04}  show that there exist balanced signature pseudo-Riemannian manifolds which are  $r$-curvature homogeneous and not locally
homogeneous for any $r$ (although the dimension of these manifolds is roughly twice $r$).  If $m := \min\{p,q\}$, then there are no known examples of
$(m+1)$-curvature homogeneous manifolds of signature $(p,q)$ which are not locally homogeneous.   These considerations have led Gilkey to conjecture \cite{GSI}  that $k_{p,q}=m+1$.

The examples in the higher signature setting above were not originally constructed for the study of curvature homogeneity, and this leads us
to a motivation for this study.  In fact, the manifolds in \cite{DG, DGS, GS04}   appeared in \cite{GIZ}, and the manifolds in
\cite{GS} appeared in \cite{GSs}--they were used  as counterexamples to the Osserman conjecture \cite{GBk, GKV} in the higher signature
setting.  As a result, the known examples have very rigid signatures.  The manifolds in \cite{DG, DGS} have balanced signature, and the
manifolds in \cite{GS} have signature $(2s,s)$ for $s \geq 1$.  It is the aim of this article to provide examples in the higher signature
setting of a more arbitrary signature.

The following is an example of a $0$-model that will be central to our discussion. 
\begin{definition} \label{theorem:def}
\rm  Let $k \geq 1$ be an integer, and choose $a, b \in \mathbb{N} \bigcup \{0\}$ so that $a+b = k$. Let $\ep_i$ be a choice of signs.  Let 
$\{U_0,...,U_k,V_0,...,V_k,S_1,...,S_k\}$
be a basis for $\mathbb{R}^{3k+2}$.  For $i = 1, \ldots, k$, we define the nonzero entries of a symmetric nondegenerate bilinear form $(\cdot,\cdot)$  and algebraic curvature tensor $R$ on the basis above as:
\begin{equation} \label{equation:0model}
(U_0, V_0) = (U_i,V_i)=1,\quad (S_i,S_i)=\ep_i,\quad\hbox{and}\quad
R(U_0, U_i, U_i, S_i) = 1\,.
\end{equation}
We define the $0$-model $\ms:= (\mathbb{R}^{3k+2}, (\cdot,\cdot), R)$.  Let $\gv$ be the structure group of this
$0$-model.  We define a \emph{normalized basis} for $V$ to be a basis that preserves the normalizations given in Equation (\ref{equation:0model}).  Thus the structure group $\gv$ can be viewed as the set of normalized bases for $\ms$.  \hfill $\qedbox$
\end{definition}
Using the same $k, a, b,$ and $\ep_i$ in Definition \ref{theorem:def}, we now define a family of pseudo-Riemannian manifolds.   
\begin{definition}  \label{theorem:manifolds} \rm
 Put coordinates
$(u_0, \ldots, u_k, v_0,
\ldots, v_k, s_1, \ldots, s_k)$ on the Euclidean space $M:=\mathbb{R}^{3k+2}$.  Let $F:=(f_1(u_1),...,f_k(u_k))$ where
$f_i(u_i)$ are a collection of smooth functions with $f_i(u_i)+1\ne0$ for all $u_i$. Define the nonzero entries of
a symmetric metric $g_F$ on the coordinate frames as follows:
$$\begin{array}  {l l}
g_F(\pu{0}, \pu{i}) = 2f_i(u_i) s_i, &\qquad g_F(\pu{i}, \pu{i}) = 
-2u_0 s_i,  \\
g_F(\pu{i}, \pv{j}) = \delta_{ij}, &\qquad g_F(\ps{i}, \ps{i}) = \ep_i\,.
\end{array}$$
Let $\mathcal{M}_F :=
(\mathbb{R}^{3k+2}, g_F).$  If we choose $a$ of the $\ep_i$ to be $-1$ and $k-a
= b$ of the $\ep_i$ to be $+1$, then this is a manifold of signature $(k+1+a,
k+1+b)$. \hfill $\qedbox$
 \end{definition}

We shall show that the manifolds $\mf$ are $0$-curvature homogeneous:
\begin{theorem}  \label{theorem:thm-1.4}
 Adopt the notation of Definition \ref{theorem:def} and of Definition \ref{theorem:manifolds}.  The manifolds $\mf$ are $0$-modeled on $\ms$.
\end{theorem}

 Define the subspaces of the model space
$V$ as follows:  
\begin{equation} \label{equation:av}
A_V :=
\{\xi\in V | R(\xi, *,*,*) = 0\} = \ker(R),  \quad A_{S,V} := A_V^\perp.
\end{equation}
 These spaces are necessarily preserved by any isomorphism of the structure group because they are defined in a basis-free
fashion. We will prove the following result involving the group of permutations $\operatorname{Sym}_k$ of $k$ objects that reflects the rigid nature of
this group:
\begin{theorem}  \label{theorem:thm-1.5}
Adopt the notation of Definition \ref{theorem:def}. If
$A$ is an isomorphism of $\ms$, then there exists a permutation
$\sigma\in\operatorname{Sym}_k$ and constants
$a_0$, $b_i$ with
$|a_0|b_i^2=1$ so that
$$\begin{array}{ll}
A U_0=a_0U_0+\Xi_0&\hbox{for some}\quad \Xi_0\in A_V,\\
A U_i=b_iU_{\sigma(i)}+\Xi_i\quad&\hbox{for some}\quad\Xi_i\in A_{S,V},\\
A S_i=\operatorname{sign}(a_0) S_{\sigma(i)}+\bar\Xi_i\quad& \hbox{for
some}\quad
\bar\Xi_i\in A_V\,.
\end{array}$$\end{theorem}

A natural question to ask is whether or not the manifolds $\mf$ are really built from
smaller dimensional manifolds with the same properties.  We recall some basic definitions relevant to this question.

\begin{definition}  \label{theorem:decompose} \rm
We say that  a $k$-model $\ms_k =(V,(\cdot,\cdot),A^0,...,A^k)$ is {\it
decomposable} if there exists a non-trivial orthogonal decomposition
$V=V_1\oplus V_2$ which induces an orthogonal decomposition $A^i=A_1^i\oplus
A_2^i$ for
$0\le i\le k$;
in this setting, we shall write $\ms=\ms^1\oplus\ms^2$
where the
$k$-model $\ms^p:=(V_p,(\cdot,\cdot)|_{V_p},A_p^0,...,A_p^k)$ for
$p = 1$ and $2$. One says
that
$\ms_k$ is {\it indecomposable} if
$\ms_k$ is not decomposable.  One says that a smooth pseudo-Riemannian manifold
$\mathcal{M}$ is {\it locally
decomposable} at a point
$P\in M$ if there exists a neighborhood $\mathcal{O}$ of $P$ so that
$(\mathcal{O},g_M)=(\mathcal{O}_1\times \mathcal{O}_2,g_1\oplus g_2)$
decomposes as a Cartesian product. We say
$\mathcal{M}$ is {\it locally
indecomposable} at $P$ if this does not happen.  \hfill $\qedbox$

It is easy to see that  if
$\ms_k(\mathcal{M},P)$ is indecomposable for some $k$, then
$\mathcal{M}$ is
locally indecomposable at $P$.   We shall show that the manifolds $\mf$ are locally indecomposable at every point in Theorem \ref{theorem:thm-1.7}:
\end{definition}
\begin{theorem}   \label{theorem:thm-1.7}
 Adopt the notation of Definition \ref{theorem:def} and of Definition \ref{theorem:manifolds}. 
 \begin{enumerate}
\item The model space $\ms$ is indecomposable.
\item The  manifolds $\mathcal{M}_F$ are locally indecomposable at every point.
\end{enumerate}
\end{theorem}
Using Theorem \ref{theorem:thm-1.5}, we can produce new isometry invariants which are not of Weyl type.  For example, in Section 5 we prove the following:
\begin{theorem}\label{theorem:extra}
Adopt the notation of Definition \ref{theorem:manifolds}.
\begin{enumerate}
\item  The following quantity is an $\ell$-model invariant:
$$
\beta_\ell  =   \sum_{j = 1}^k  \frac{f_j^{(\ell +1)} (1+f_j')^{\ell-1}
}{\left[f_j^{(2)}\right]^{\ell}}\,.
$$
\item  If the manifold $\mathcal{M}_F$ is $\ell$-curvature homogeneous, then
$\beta_p$ is constant for all $p = 1, 2, \ldots, \ell$.

\item  If $\mathcal{M}_F$ is locally homogeneous, then $\beta_\ell$ is constant
for all $\ell$.

\end{enumerate}
\end{theorem}
Using this theorem and a similarly defined $\ell$-model invariant (see Theorem \ref{theorem:gammas}), it is possible to prove:
\begin{theorem}  \label{theorem:thm-1.8}
Suppose $f_i'(u_i)+1 \neq 0$ for $1 \leq i\leq k$.  If $f_i''(u_i)\neq 0,$ then $\mf$ is not $2$-curvature homogeneous.
\end{theorem}

The following is a brief outline of the paper.  We will compute the entries of tensors $R$ and $\nabla R$, and prove Theorem \ref{theorem:thm-1.4} in Section \ref{section:two}.  In Section \ref{section:three} we study the structure group $\gv$ and establish Theorem \ref{theorem:thm-1.5}.  We study the notion of indecomposability in Section \ref{section:four} and prove Theorem \ref{theorem:thm-1.7}.  In Section \ref{section:five} we conclude the paper by establishing Theorems \ref{theorem:extra} and  \ref{theorem:thm-1.8}.

\section{Curvature Homogeneity}  \label{section:two}

We begin this section with a calculation of the Christoffel symbols of the Levi-Civita connection of the manifolds $\mf$.  
\begin{lemma}   \label{theorem:calc}
Let $\pu i, \ps i$ and $\pv i$ be coordinate vector fields on $\mathcal{M}_F$.

 \begin{enumerate}

\item  The nonzero covariant derivatives of the
coordinate vector fields are

$$\begin{array}{r c l}
\nabla_{\pu{0}} \pu{i} = \nabla_{\pu{i}} \pu{0} & = &-s_i\pv{i}
-f_i(u_i)
\ep_i
\ps{i}, \\
\nabla_{\pu{i}} \pu{i} & = & (2f'_i(u_i)+1)s_i\pv{0} +
u_0 \ep_i\ps{i},
\\
\nabla_{\pu{0}} \ps{i} = \nabla_{\ps{i}} \pu{0} & = & f_i(u_i) \pv{i},  \\
\nabla_{\pu{i}} \ps{i} = \nabla_{\ps{i}} \pu{i} & = & f_i(u_i) \pv{0}-
u_0\pv{i}\,.
\end{array}$$

\item  The only nonzero entries of the Riemannian curvature tensor $R$ (up to
the usual
$\mathbb{Z}_2$ symmetries) are

\begin{enumerate}
\item  $R_0(i) := R(\pu{0}, \pu{i}, \pu{i}, \pu{0}) =
f_i(u_i)^2\ep_i,$ and  
\item  $R_s(i):= R(\pu{0}, \pu{i}, \pu{i}, \ps{i}) = f_i'(u_i)+1.$
\end{enumerate}

\item  The only nonzero entries of the covariant derivative tensor $\nabla
R$  (up to the usual symmetries)  are:

\begin{enumerate}
\item  $\nabla R(\pu{0}, \pu{i}, \pu{i}, \pu{0} ; \pu{i}) = 2f_i(u_i)
\ep_i(2f_i'(u_i) + 1)$
\item  $\nabla R(\pu{0}, \pu{i}, \pu{i}, \ps{i} ; \pu{i}) = f_i''(u_i)$
\end{enumerate}

\item  The following assertions are equivalent:
\begin{enumerate}
\item For each $i$ with $1 \leq i \leq s$, either $f_i(u_i) = 0$ or
$f_i^\prime(u_i) = -\frac{1}{2}$.
\item $\mathcal{M}_F$ is a local symmetric space.
\end{enumerate}

\end{enumerate}
\end{lemma}
\begin{proof} 

We compute the nonzero components of the covariant  derivatives of the
coordinate vector fields, the curvature tensor $R$ and its
covariant derivative $\nabla R$.  Note that
$g(\pu{j},
\ps{i}) = g(\pv{j},\ps{i}) = 0$ and $g(\ps{i}, \ps{i}) = \ep_i$ is 
constant.  So if $X$ and $Y$ are any coordinate vector fields, we have
$$g(\nabla_{\ps{i}} X, Y) = g(\nabla_{X}\ps{i}, Y) = -g(\nabla_XY, \ps{i}) =
\frac{1}{2}(\ps{i}g(X, Y))\,.$$  We let the index $i$ range from $1$ to $k$.
\begin{equation*}
\begin{array}{r c l}
g(\nabla_{\pu{0}}\pu{i}, \pu{i}) & = &\frac{1}{2}\pu{0}g(\pu{i}, \pu{i}) \\
  & = & \frac{1}{2}(-2s_i) = -s_i,   \\
g(\nabla_{\pu{0}}\pu{i}, \ps{i}) & = & \frac{1}{2}( \pu{0}g(\pu{i}, \ps{i}) + \pu{i}
(\pu{0}, \ps{i}) - \ps{i} g(\pu{0}, \pu{i}))  \\
  & = & \frac{1}{2}(2f_i) = f_i,     \\
g(\nabla_{\pu{i}}\pu{i}, \pu{0}) & = & \frac{1}{2}( 2\pu{i}g(\pu{i}, \pu{0})  -
\pu{0} g(\pu{i}, \pu{i}))  \\
  & = &\frac{1}{2}(2\cdot2f_i's_i - (-2s_i) = s_i(2f_i'+1),     \\
 g(\nabla_{\pu{i}} \pu{i}, \ps{i}) & = &
-\frac{1}{2}(\ps{i}g(\pu{i},\pu{i})=u_0, \\ 
g(\nabla_{\pu{0}} \ps{i}, \pu{i}) & = & \frac{1}{2}(\ps{i} g(\pu{0}, \pu{i}))
= f_i,
\\ g(\nabla_{\pu{i}} \ps{i}, \pu{0})   & = & \frac{1}{2}(\ps{i} g(\pu{i},
\pu{0})) = f_i, \\
 g(\nabla_{\pu{i}} \ps{i}, \pu{i}) & = & \frac{1}{2}(\ps{i}g(\pu{i}, \pu{i})) =
-u_0\,.
\end{array}\end{equation*} 
We may then use this computation to
see that: 
\begin{equation*} \begin{array}{l}
 R(\pu{0}, \pu{i}) \pu{i}=(\nabla_{\pu{0}}\nabla_{\pu{i}} -\nabla_{\pu{i}}
\nabla_{\pu{0}} ) \pu{i}  \\
\quad= \nabla_{\pu{0}}[(2f_i'+1)s_i \pv{0} + u_0 \ep_i \ps{i}] -
\nabla_{\pu{i}}[-s_i \pv{i} - f_i \ep_i \ps{i}] \\
\quad= \ep_i \ps{i} + u_0 \ep_i\nabla_{\pu{0}}\ps{i} + f_i' \ep_i \ps{i} + f_i \ep_i
\nabla_{\pu{i}} \ps{i}  \\
\quad= (1+f_i')\ep_i \ps{i} + f_i^2 \ep_i \pv{0}\,.\end{array}
\end{equation*}
The covariant derivative of $R$ is given by:
\begin{equation*}\begin{array}{l}
\nabla R(\pu{0}, \pu{i}, \pu{i}, \pu{0}; \pu{i})\\\quad = \pu{i}(f_i^2 \ep_i) -
2R(\nabla_{\pu{i}} \pu{0}, \pu{i}, \pu{i}, \pu{0})  - 2R(\pu{0},
\nabla_{\pu{i}}\pu{i}, \pu{i}, \pu{0}) \\
\quad=2f_if_i'\ep_i + 2f_i\ep_i(f_i'+1) = 2f_i\ep_i(2f_i'+1),  \\
\nabla R (\pu{0}, \pu{i}, \pu{i}, \ps{i}; \pu{i})\\\quad=
\pu{i}(f_i'+1) - R(\nabla_{\pu{i}}\pu{0}, \pu{i}, \pu{i}, \ps{i}) - R(\pu{0},
\nabla_{\pu{i}}\pu{i}, \pu{i}, \ps{i})  \\
\qquad - R(\pu{0}, \pu{i}, \nabla_{\pu{i}}\pu{i}, \ps{i}) - R(\pu{0}, \pu{i}, \pu{i},
\nabla_{\pu{i}}\ps{i})  \\
\quad= f_i''\,.
\end{array}
\end{equation*}
The Lemma now follows.
\end{proof}
We establish Theorem \ref{theorem:thm-1.4} after a brief remark.
\begin{remark} \label{remark:gpwm}  \rm
Let the index $\mu$ range from $1$ to $k$, and let the index $\nu$ range from $0$ to $k$.  If we relabel the coordinates $x_\nu = u_\nu$, $x_{k+\mu} = s_\mu$, and $x_{2k+1+\nu} = v_\nu$, the above calculations show that $\nabla_{\partial_{x_i}} \partial_{x_j} = \sum_{k > \max\{i,j\}} \Gamma_{ij}{}^k(x_0, \ldots, x_{k-1}) \partial_{x_k}$.  Thus by definition, $\mf$ is a family of \emph{generalized plane wave manifolds}.  By the results of Gilkey and Nik\v cevi\'c \cite{G-S}, we conclude that members of the family $\mf$ are Ricci-flat, complete, $\exp:T_PM \to M$ is a diffeomorphism for all $P$, and all Weyl scalar invariants vanish.  We will see in Section \ref{section:five} that there are members of the family $\mf$ which are not locally homogeneous.    This is not possible in the Riemannian setting as Pr\"ufer, Tricerri, and Vanhecke \cite{PTV} showed that if all local scalar Weyl invariants up to order $\frac12n(n-1)$ are constant on a Riemannian manifold $(N,h)$ of dimension $n$, then $(N,h)$ is locally homogeneous and determined up to local isometry by these invariants.  \hfill $\qedbox$
\end{remark}
\emph{Proof of Theorem \ref{theorem:thm-1.4}. }  To show that $\mf$ are $0$-modeled on $\ms$, we will produce a normalized basis for $(T_PM, g|_P, R|_P)$ for any $P \in M$ (see Definition \ref{theorem:def}).  We have that $f_i(u_i)+1\ne0$ for $1\le i\le k$.  We set
\begin{equation*}
\begin{array}{lll}
U_0:=\partial_{u_0}+\sum_ja_j\partial_{s_j},\quad
&U_i:=b_i\partial_{u_i}+\beta_i\partial_{v_0}+\tilde\beta_i\partial_{v_i},\\
S_i:= \kappa_i\partial_{s_i}+\gamma_i\partial_{v_i},&
V_0:=\partial_{v_0},\\
V_i=b_i^{-1}\partial_{v_i},\end{array}
\end{equation*}
where $b_i$, $\beta_i$, $\tilde \beta_i$, $\kappa_i$, and
$\gamma_i$ will be specified presently. The potentially non-zero curvatures
are then:
$$
\begin{array}{l}
R(U_0,U_i,U_i,U_0)=b_i^2\{f_i(u_i)^2\varepsilon_i+2a_i(f_i^\prime(u_i)+1)\},\\
R(U_0,U_i,U_i,S_i)=b_i^2(f_i^\prime(u_i)+1)\varepsilon_i\kappa_i\,.
\end{array}
$$
To ensure that $R(U_0,U_i,U_i,U_0)=0$ and $R(U_0,U_i,U_i,S_i)=+1$, we set
\begin{equation*}
\begin{array}{l}
a_i:=-\frac{f_i(u_i)^2\varepsilon_i}{2(f_i^\prime(u_i)+1)},\\
  \kappa_i:=\varepsilon_i\operatorname{sign}(f_i^\prime(u_i)+1),
\vphantom{\vrule height 11pt}\\
  b_i:=|f_i^\prime(u_i)+1|^{-1/2}\,.
\vphantom{\vrule height 11pt}\end{array}\end{equation*}
The potentially non-zero inner products are
$$\begin{array}{ll}
(U_0,V_0)=1,&(U_0,S_i)=\kappa_ia_i+\gamma_i,\\
(U_0,U_i)=b_ig_F(\partial_{u_0},\partial_{u_i})+\beta_i,&(S_i,S_i)=1,\\
(U_i,U_i)=b_i^2g_F(\partial_{u_i},\partial_{u_i})+2b_i\tilde\beta_i,
\quad&(U_i,V_i)=1
\,.\vphantom{\vrule height 11pt}
\end{array}
$$
We complete the proof by setting:
$$
\begin{array}{l}
\gamma_i:=-\kappa_ia_i,\quad\beta_i:=-b_ig_F(\partial_{u_0},\partial_{u_i}),\\
\tilde\beta_i:=-\textstyle\frac12b_ig_F(\partial_{u_i},\partial_{u_i})\,.
\end{array}
$$
\hfill $\qedbox$

It will be convenient to compute several values of the curvature tensor and its
covariant derivatives on a normalized basis, see Theorems \ref{theorem:independent} and \ref{theorem:gammas}.  We list these quantities
below for future reference.

\begin{lemma} \label{theorem:nabla-on-norm}
Adopt the notation of Definition \ref{theorem:def} and Definition \ref{theorem:manifolds}.  Suppose that  $\{U_0, \ldots, U_k, V_0, \ldots, V_k, S_1, \ldots, S_k\}$ is the normalized basis found in the previous theorem.
\begin{enumerate}
\item $ \nabla R(U_0, U_i, U_i, U_0; U_i) =  \frac{f_i\ep_i}{(f_i'+1)^{5/2}}\left[ 2(2f_i'+1)(f_i'+1) - f_if_i''
\right]$.
\item  $\nabla R(U_0, U_i, U_i, S_i; U_i)=  \frac{f_i''\kappa_i}{|f_i'+1|^{3/2}}$.
\item  $\nabla^{\ell} R(U_0, U_i, U_i, S_i; U_i,\ldots, U_i) =  \kappa_if_i^{(\ell+1)}|f_i'+1|^{-\frac{2+\ell}{2}} $
\item  $\nabla^2R(U_0, U_i, U_i, U_0; U_i, U_i) = \frac{\varepsilon_i}{(f_i^\prime+1)^2}
\left(
4(f_i^\prime)^2 +2f_i^\prime + 6f_if_i^{\prime\prime} -
\frac{(f_i)^2f_i^{\prime\prime\prime}}{f_i^\prime+1}
\right)\,.$
\end{enumerate}
\end{lemma}

\begin{proof}  We use the normalized basis
found in the proof of Theorem \ref{theorem:thm-1.4} and the calculations of Lemma \ref{theorem:calc} to compute these directly--the calculations are omitted.\end{proof}

\section{The Structure Group $\gv$}  \label{section:three}

In this section we study the structure group $\gv$.  For convenience, we establish notation as follows for the normalized bases $\BB$ and $\tilde\BB$: 
$$\begin{array}{c}
\BB = \{U_0, \ldots, U_k, V_0, \ldots, V_k, S_1,\ldots S_k\}, \\
\tilde\BB = \{\tilde U_0, \ldots, \tilde U_k, \tilde V_0, \ldots, \tilde V_k,
\tilde S_1,\ldots \tilde S_k\}.
\end{array}
$$ 
 We adopt the notation of Equation (\ref{equation:av}).  For any normalized basis $\BB$, one has 
\begin{eqnarray*}
A_{V} &=& \CDspan\{V_0, \ldots, V_k\}, \hbox{ and } \\
A_{S,V}&=& \CDspan\{S_1, \ldots, S_k, V_0, \ldots, V_k\}\,.
\end{eqnarray*}
  Let $\operatorname{Sym}_k$ be the group of permutations of the numbers $\{1,
\ldots, k\}$.  

\emph{Proof of Theorem \ref{theorem:thm-1.5}.}
Note $A S_i\in A_{S,V}$. We expand:
\begin{equation}\label{eqn-3.a}
\begin{array}{r c c c l}
A U_0&=&a_0U_0 &+&\sum_j(b_{0j}U_j+d_{0j}S_j)+A_V,\\
A S_i&=& & & \sum_jf_{ij}S_j+ A_V, \\
A U_i&=&a_iU_0&+&\sum_jb_{ij}U_j+A_{S,V}\,.
\end{array}\end{equation}

For any $\xi_1,\xi_2\in V$, we have that:
\begin{equation}\label{eqn-3.b}
0=R(\xi_1,U_0,U_0,\xi_2)=R(A\xi_1,A U_0,A U_0,A\xi_2)\,.
\end{equation}
Choose $\xi_i$ so $A\xi_1=U_0$ and $A\xi_2=S_j$. We then have
$$
0=R(U_0,A U_0,A U_0,S_j)=b_{0j}^2\,.
$$
Consequently $b_{0j}=0$.
We have $A\cdot A_V=A_V$. As
$1=( U_0,V_0)=(A U_0,A V_0)$, there exists 
$v\in A_V$ so $(A  U_0,v)\neq 0$. Since $A
U_0=a_0U_0+A_{S,V}$, we conclude
$a_0\ne0$. Choosing  $A\xi_1=A\xi_2=U_i$ in Equation
(\ref{eqn-3.b}) we have:
$$0=R(U_i,A U_0,A U_0,U_i)=2a_0d_{0j}\,.$$
Since $a_0\ne0$, $d_{0j}=0$.  Display (\ref{eqn-3.a}) becomes
\begin{eqnarray*}
&&A U_0=a_0U_0+A_V,\quad
  A S_i=\sum_jf_{ij}S_j+A_V,\\
&&A U_i=a_iU_0+\sum_jb_{ij}U_j+A_{S,V}\,.
\end{eqnarray*}

Since $A V_i\in A_V$, the matrix $[b_{ij}]$ is invertible. 
 Suppose the matrix element $b_{ij}\ne0$. Choose $\xi_1$ so $A\xi_1=S_j$.  Since
$k \geq 2$, we may choose positive induces $\ell\ne i$, then
 $$0=R(U_0,U_i,U_\ell,\xi_1)=R(A U_0,A U_i,A
U_\ell,A\xi_1)=a_0b_{ij}b_{\ell j}\,.$$ Thus if $b_{ij}\ne0$, $b_{\ell j}=0$
for
$i\ne \ell$. So in the matrix $b_{ij}$, each column has at most one non-zero
entry. Since
$b_{ij}$ is invertible, each column has exactly one non-zero entry. So
one has:
\begin{eqnarray*}
&&A U_0=a_0U_0+ A_V,\quad
  A S_i=\sum_jf_{ij}S_j+ A_V,\\
&&A U_i=a_iU_0+b_iU_{\sigma(i)}+ A_{S,V}\,.
\end{eqnarray*}
The relation $\delta_{ij}=R(A U_0,A U_i,A U_i,A S_j)$ shows
$f_{ij}=0$ for $j\ne\sigma(i)$. Since $A S_j$ is a unit vector, this
coefficient is $\pm1$. Thus
$$A U_0=a_0U_0+A_V,\quad
  A S_i=\pm S_{\sigma(i)}+A_V,\quad
A U_i=a_iU_0+b_iU_{\sigma(i)}+ A_{S,V}\,.
$$
Since $1=R(A U_0,A U_i,A U_i,A S_i)$, we have  $\pm
b_i^2a_0=1$. Finally, since $k\ge2$ and since $0=R(A U_i,A U_j,A
U_j,A S_i)$, we have $a_ib_j=0$ and hence $a_i=0$. The relation
$|a_0|b_i^2=1$ and
$A S_i=\operatorname{sign}(a_0)S_{\sigma{(i)}}$ now follow.  This
establishes the theorem.\hfill $\qedbox$
\begin{remark}  \rm
Theorem \ref{theorem:thm-1.5} does not apply when $k = 1,$ although similar
statement is true in that case:  If $A$ is an isomorphism of
$\ms$ then
$$\begin{array}{ll}
A U_0=a_0U_0+\Xi_0&\hbox{for some}\quad \Xi_0\in A_V,\\
A U_1=a_1U_0 + b_1U_1+\Xi_1\quad&\hbox{for some}\quad\Xi_1\in
A_{S,V},\\
A S_1=\operatorname{sign}(a_0) S_1+\bar\Xi_1\quad& \hbox{for
some}\quad
\bar\Xi_1\in A_V\,.
\end{array}$$
Notice the extra freedom in choosing $a_1$.   Since $\operatorname{Sym}_1$ is the trivial group, the symmetric group action is not so evident as when $k \geq 2$. \hfill $\qedbox$
\end{remark}

The crucial part of the previous result is that any change of basis will
permute the interesting information, single out the vector $U_0$ and
$A\cdot A_{S,V} \subseteq A_{S,V}$.  This will be important when
defining  invariants in the next section.  The extra information one has when $k
= 1$ will not create any ambiguity in the development of any of our invariants.

\section{Indecomposability}  \label{section:four}

Since $\mathbb{R}^{3k+2}$ is contractible, any real vector bundle over
$\mathbb{R}^{3k+2}$ is trivial, in particular, the tangent bundle is trivial.  With the
added structure of a metric and a curvature tensor, however, more information
is available.  

 A natural question to ask is if these
manifolds are really products of manifolds of smaller dimension.  More
specifically, is
$\mathbb{R}^{3k+2} = M_1 \times M_2$ and $g_F = g_{M_1}
\oplus g_{M_2}$?  If this were the case, then  $T\mathbb{R}^{3k+2} = TM_1
\oplus TM_2$, and one has that the curvature tensor $R_M = R_{M_1} \oplus
R_{M_2}$.  This is a more algebraic notion of indecomposability which we
briefly study.  The motivation comes from the main result in \cite{TV}:  any family of Riemannian manifolds $0$-modeled on an irreducible symmetric space are homogeneous (in fact, symmetric).  In the pseudo-Riemannian setting, the notion of irreducibility seems more elusive, and although we do not show that the $0$-model $\ms$ is irreducible, we prove the weaker Theorem \ref{theorem:thm-1.7}.  Although, the main step of the result in \cite{TV} is to use the hypothesis to establish that the manifolds in question are Einstein.  We recall Remark \ref{remark:gpwm}:  the manifolds $\mf$ are not only Einstein, but Ricci-flat.  Thus this family of manifolds provide interesting insight into the distinction between Riemannian and pseudo-Riemannian manifolds.

 Recall the notation established in Definitions \ref{theorem:def} and \ref{theorem:manifolds}.  We show in this section that the
manifolds
$\mathcal{M}_F$ are locally indecomposable at every point, and thus locally
$\mathcal{M}_F$ is not the direct product of smaller dimensional manifolds,
answering the above question in the negative.  

We fix a normalized basis
$\BB$ for this section. Using the subspace $A_V$ defined in the introduction, denote $V/A_V = B_{U,S}$, and
$\pi:V \rightarrow B_{U,S}$ the projection.  A basis for
$B_{U,S}$ is the image of $U_0, \ldots, U_k, S_1, \ldots, S_k$ under $\pi$. 
Write $\bar U_i = \pi U_i$, similarly for the other vectors.  Since 
$A_V\subset\ker(R)$, we have a well-defined algebraic curvature tensor $\bar
R$ defined on
$B_{U,S}$, characterized by the relation
$\pi^* 
\bar R = R$.  We have the same relations for
$\bar R$ on the image of the normalized basis as we do for $R$ on the 
original normalized basis for $V$, although of course the projection of such
a basis to $B_{U,S}$ is no longer linearly independent.   We recall that on $V$,
we have the relations
 $$( U_i,V_i)=\delta_{ij},\quad(S_i,S_i)=\varepsilon_i,\quad R(U_0,U_i,U_i,S_i)=1\,.$$

\begin{lemma}  \label{theorem:D-model}
The weak $0$-model $(B_{U,S}, \bar R)$ is indecomposable for $k \geq 1$.
\end{lemma} 
\begin{proof}  We assume to
the contrary there exists a non-trivial decomposition of the model space $(W,R)=(\bar W_1\oplus\bar W_2,R_1\oplus R_2)$ and argue for a
contradiction. We begin by expressing
$\bar U_0 =
\xi_1 + \xi_2$, for $\xi_i \in W_i$.  

\emph{Case I.}  One of $\xi_i$ is 0   (suppose without loss of generality
that $\xi_2 = 0$).  This means that we can write $\bar U_0 \in \bar W_1$.  Let
$0 \neq \eta \in \bar W_2.$  Consequently, we may express $\eta = \gamma_0 \bar
U_0 +
\sum_{j = 1}^k
\gamma_j
\bar U_j + \gamma_j' \bar S_j$.  Then for $i > 0$, 
$$
\begin{array}{r c l}
\bar R(\bar U_0, \bar U_i, \bar U_i, \eta) & = & \gamma_i' = 0,\hbox{ and }  \\
\bar R(\bar U_0, \bar U_i, \eta, \bar S_i) & = & \gamma_i = 0\,.
\end{array}
$$
So $\eta = \gamma_0 \bar U_0$, and $\eta \neq 0$ means that $\eta \in W_2$ and
$U_0 \in \bar W_1$ are not linearly independent, and so 
$W_1\cap W_2\ne\{0\}$.  This contradiction permits us to eliminate this
case from consideration. 

\emph{Case II.}  $\bar U_0 = \xi_1 + \xi_2$ and both $\xi_i \neq 0$.  We
express these vectors as 
$$
\begin{array}{r c l}
\xi_1 & = & \alpha_0 \bar U_0 + \sum_j \alpha_j \bar U_j + \alpha_j' \bar S_j,
\\
\xi_2 & = & \beta_0 \bar U_0 + \sum_j \beta_j \bar U_j + \beta_j' \bar S_j\,.
\end{array}
$$
Since $\xi_1 + \xi_2 = \bar U_0$, we must have $\alpha_0 + \beta_0
= 1$, $\alpha_j + \beta_j = \alpha_j' + \beta_j' = 0$.  For $j = 1, 2$ and $i
= 1, \ldots, k$, we compute
$\bar R (\bar U_0, \xi_j, \xi_j, \bar S_i)$ in two ways.  First, we could have
only the $\bar U_i$ coefficients of $\xi_j$, so $\bar R (\bar U_0, \xi_j,
\xi_j,
\bar S_i) = \alpha_i^2$ ($j = 1$)  or $\beta_i^2$  ($j = 2$).  On the other
hand (for $j = 1$), 
$$
\begin{array}{r c l}
\bar R (\bar U_0, \xi_1, \xi_1, \bar S_i) & = & \bar R (\xi_1 + \xi_2, \xi_1,
\xi_1, \bar S_i)  \\
 & = &   \bar R (\xi_1, \xi_1,
\xi_1, \bar S_i)  + \bar R (\xi_2, \xi_1,
\xi_1, \bar S_i) \\
  & = &  0.
\end{array}
$$
Similarly for $j = 2$.  Thus $\alpha_i = \beta_i = 0$ for all $i$.  

Now we go to work on the other coefficients.  Since $\alpha_0 + \beta_0 = 1$,
at least one of these must be nonzero.  Suppose without loss of generality
that $\alpha_0 \neq 0$.    Compute $0=\bar R(\xi_1, \bar U_j, \bar U_j,
\xi_2) =
\alpha_0 \beta_j' + \beta_0 \alpha_j'$.  Since $\alpha_0 \neq 0$, we can solve
for $\beta_j' = \frac{-\beta_0 \alpha_j'}{\alpha_0}$.  Imposing the condition
$\alpha_j' + \beta_j' = 0$ gives us $\alpha_j'(\alpha_0 - \beta_0)= 0$ for
all $j = 1, 2, \ldots, k$.  These equations could be solved by having either
$\alpha_j' = 0$ for all $j$ or 
$\alpha_0 = \beta_0 $.  

\emph{Case II.a.} Suppose we have $\alpha_j' = 0$ for all $j$.  Then we 
again impose the condition $\alpha_j' + \beta_j' = 0$ to see that $\beta_j' =
0$ for all $j$ as well.   This gives us $\xi_1 = \alpha_0 \bar U_0$ and
$\xi_2 = \beta_0 \bar U_0$, and at this point there are several contradictions: 
by assumption, both $\xi_i$ are nonzero, and we have
$\xi_1 =
\lambda \xi_2$, not linearly independent, but living in different subspaces. 
This is false. 

\emph{Case II.b.}  Suppose $\alpha_0 = \beta_0$.  Then $\alpha_0 + \beta_0 = 1$
implies $\alpha_0 + \beta_0 = \frac{1}{2}$.  Unfortunately, we must go into
further cases and consider where another vector lives.  The analysis of this
new vector is similar to the previous technique. Since $k
\geq 1$, there exists a $\bar U_1 \in B_{U,S}$, and we proceed by studying 
$\bar U_1$.  Write $\bar U_1 = \eta_1 + \eta_2,$
and
$\eta_i \in \bar W_i$.

\emph{Case II.b.i.}  One of $\eta_i = 0$.  Without loss of generality, assume
$\eta_2 = 0$.  Then $\bar U_1 \in \bar W_1$.  Then $\bar R (\xi_2, \bar U_1,
\bar U_1, \bar S_1) = \frac{1}{2}$, but since $\xi_2 \in \bar W_2$ and $\bar
U_1 \in \bar W_1$, we must have $\bar R (\xi_2, \bar U_1,
\bar U_1, \bar S_1) = 0$ which gives us a contradiction.

\emph{Case II.b.ii.}  Both $\eta_i \neq 0$.  We write $\eta_i = a_i \bar U_1 +
v_i$ for $v_i \in \bar W_i$.  Then $a_1 + a_2 = 1$ and hence both $a_i$ cannot
be 0 simultaneously.  We compute 
$$
\begin{array}{r c l}
\bar R(\xi_2, \bar U_1, \eta_1, \bar S_1) & = & \frac{1}{2}a_1 =  0,  \\
\bar R(\xi_1, \bar U_1, \eta_2, \bar S_1) & = & \frac{1}{2}a_2 =  0\,.
\end{array}
$$
This yields
a contradiction; this final contradiction completes the proof.\end{proof}  
\noindent\emph{Proof of Theorem \ref{theorem:thm-1.7}.}   We have shown in Lemma \ref{theorem:D-model} that the weak model space $B_{U,S}$ is indecomposable.  In addition, $\ker R = \CDspan\{V_0, \ldots, V_k\}$ is a totally isotropic subspace.  Thus according to \cite{GBk2}, the model space $\ms$ is indecomposable.

 We now prove Assertion (2).  We have shown that $\ms$ is
a $0$-model for the tangent space $T_PM$ at any point $P \in M$.  Such a
decomposition of $T_PM$ would induce a decomposition of the
$0$-model $\ms$.  But $\ms$ is indecomposable
by Assertion (1), and no such decomposition of the tangent bundle is
possible. \hfill $\qedbox$

\section{Isometry Invariants and Local Homogeneity}  \label{section:five}

Since all Weyl scalar invariants vanish (see Remark \ref{remark:gpwm}) we use the determination of the structure group $\gv$ given in Theorem \ref{theorem:thm-1.5} to define new isometry
invariants.  We build invariants involving normalized bases and only the tensors
$\nabla R,
\ldots,
\nabla^{\ell} R$; these are so-called
$\ell$-model invariants.  This will aid us in studying the question of $\ell$-curvature homogeneity for $\ell \geq 2$ for the manifolds
$\mf$.   We will need a technical lemma describing the behavior of the higher covariant derivatives on a normalized basis.  

\begin{lemma}   \label{theorem:derivatives}
For the manifolds defined above, the following assertions hold. 
Let $\ell \geq 1$ and $i = 1, 2, \ldots , k$.

\begin{enumerate}

\item $ \nabla^\ell R(\pu{0}, \pu{i}, \pu{i}, \ps{i};
\pu{i},
\ldots,
\pu{i}) = f_i^{(\ell + 1)}(u_i).
$

\item  $\nabla^\ell R(\pu{0}, \pu{i}, \pu{i}, \pu{0}; \pu{i}, \ldots \pu{i})$
is a function of $u_i$, expressible as an algebraic combination of the
derivatives of
$f_i$.

\item  $\nabla^\ell R(*,*,*,*;*, \ldots, *, \ps{i}) = 0$.

\item  $\nabla^\ell R(*,*,*,*;*, \ldots, *, \pu{0}) = 0$.

\item  The only possible nonzero entries of the covariant derivatives of $R$
on any normalized basis are 
$$ \nabla^\ell R(U_0, U_i, U_i, S_i;
U_i,
\ldots,
U_i) \quad \hbox{and} \quad \nabla^\ell R(U_0, U_i, U_i, U_0; U_i, \ldots
U_i).$$

\end{enumerate}

\end{lemma}

\begin{proof}  Assertions 1 and 2 follow from Lemma \ref{theorem:calc}, Assertion 3.  Note
that in these terms, both are functions of only the $u_i$.  Hence to uncover
any other nonzero terms of the higher covariant derivatives other than those
ending in only $\pu{i}$, we must look to our calculation of $\nabla$ on the
coordinate frames (see Lemma \ref{theorem:calc}, Assertion 1).   Assertion 3 is now obvious,
and since $\nabla_{\pu{0}} \pu{0}=0$, we see Assertion 4 follows as well. 
As we may only build higher covariant derivatives
from $\pu{i}$ with those relations in Assertion 3 of Lemma 1.1,
and that any change of normalized basis will  permute the same positive $U_*$
and $S_*$ induces, the only nonzero higher covariant derivatives on any
normalized basis are only  those listed.\end{proof}

 Let  $\BB = \{U_0, \ldots, U_k, V_0,
\ldots V_k, S_1, \ldots, S_k\}$ be the normalized basis found in Theorem 2.1.   We   define below the
functions
$(\beta_\ell)_{\BB}$ for $\ell \geq 2$, which a priori depends on the
choice of normalized basis.  Assume for now that all denominators are nonzero. 
Define 
\begin{equation*}
(\beta_\ell)_\BB:= \sum_{j= 0}^k  \frac{\nabla^\ell R(U_0, U_j, U_j,
S_j; U_j,\ldots, U_j)}{\left( 
\nabla R(U_0, U_j, U_j, S_j; U_j)
\right)^\ell}\,.
\end{equation*}

\begin{lemma}   \label{theorem:independent}
Adopt the notation of Definitions \ref{theorem:def} and \ref{theorem:manifolds}.  If  $f_i'' \neq 0$ and $\ell \geq 2$, then  $(\beta_\ell)_\BB$ is independent of the
normalized basis chosen.
\end{lemma}

\begin{remark} \rm
The hypothesis $f_i' +1 \neq 0$ is required for a normalized basis to exist. 
The condition that $f_i'' \neq 0$ is required for the invariants $\beta_\ell$ to
exist at all, as we divide by the quantity $f_i''$ in the definition of
$\beta_\ell$.  These two hypothesis are needed only for these reasons; i.e., we
need everything to ``make sense''.  Later, we remove the restriction $f_i''
\neq 0$ in the definition of another invariant (see Theorem
\ref{theorem:gammas}).  \hfill $\qedbox$
\end{remark}

\emph{Proof of Lemma \ref{theorem:independent}.}
Let $\tilde \BB$ be another normalized basis, and
$\sigma \in \operatorname{Sym}_k$ be the corresponding permutation of the
induces found in Theorem \ref{theorem:thm-1.5}.  By Lemma
\ref{theorem:derivatives}, we know how a normalized change of basis effects the
entries of the higher covariant derivatives.  Essentially, the only change of
basis possible is a permutation of the $U_*$ and
$S_*$ basis vectors with a (nonzero) scaling factor.  So,

\begin{eqnarray*}
&&\nabla^\ell R(\tilde U_0, \tilde U_j, \tilde U_j,
\tilde S_j;
\tilde U_j,\ldots,  \tilde U_j)\\& =&
\left(\frac{\pm 1}{\sqrt{|a_0|}}\right)^\ell \nabla^\ell R(U_0, U_{\sigma(j)},
U_{\sigma(j)}, S_{\sigma(j)}; U_{\sigma(j)}, \ldots,  U_{\sigma(j)}), 
\end{eqnarray*} 
and
\begin{eqnarray*}
&&(\nabla
R(\tilde U_0, \tilde U_j, \tilde U_j, \tilde S_j; \tilde U_j))^\ell \\
&=&
\left(\frac{\pm 1}{\sqrt{|a_0|}}\right)^\ell\nabla R(U_0, U_{\sigma(j)},
U_{\sigma(j)}, S_{\sigma(j)}; U_{\sigma(j)})^\ell\,.
\end{eqnarray*}
The permutation $\sigma$ is a bijection of a finite set of induces, and
so if we put $$I = \{\sigma^{-1} (1), \ldots, \sigma^{-1} (k)  \} = \{ \ell_1,
\ldots, \ell_k\},$$ we get the rearranged (but equal) sum 

$$
\begin{array}{r c l}
(\beta_\ell)_{\tilde \BB} & = & 
\sum_{j = 1}^k  \frac{\nabla^\ell R(\tilde U_0, \tilde U_{\ell_j},
\tilde U_{\ell_j},
\tilde S_{\ell_j}; \tilde U_{\ell_j}, \ldots, \tilde U_{\ell_j})}{\left( 
\nabla R(\tilde U_0, \tilde U_{\ell_j}, \tilde U_{\ell_j}, \tilde S_{\ell_j};
\tilde U_{\ell_j})
\right)^\ell}
 \\
 & = & 
\sum_{j= 0}^k  \frac{\nabla^\ell R(U_0, U_j, U_j, S_j; U_j,
U_j)}{\left( 
\nabla R(U_0, U_j, U_j, S_j; U_j)
\right)^\ell}
  \\
  & = & (\beta_\ell)_{\BB}\,.
\end{array}
$$
Hence $(\beta_\ell)_{\BB} =  (\beta_\ell)_{\tilde \BB} =
\beta_\ell$ is independent of the basis chosen, and is an invariant of the
manifolds
$\mathcal{M}_F$.    \hfill $\qedbox$

\emph{Proof of Theorem \ref{theorem:extra}.}
 Evaluating these tensors on
a normalized basis and using Theorem \ref{theorem:derivatives} and Lemma \ref{theorem:nabla-on-norm} establishes the first assertion of Theorem \ref{theorem:extra}.

If $\mathcal{M}_F$ were $\ell$-curvature homogeneous, then there exists a $p$-model for every $p = 0, 1, \ldots, \ell$, along with a normalized basis for
$T_PM$ so that the metric, and curvature entries up to order $\ell$ are
constant.  Since $\beta_p$ is built from these entries, $\beta_p$  
must be constant for all $p = 0, \ldots, \ell$.   This establishes
Assertion 2 of Theorem \ref{theorem:extra}.

If $\mathcal{M}_F$ is locally homogeneous, then it is $\ell$-curvature
homogeneous for all $\ell$.  Applying the previous assertion shows that $\beta_\ell$ 
has to be constant for all $\ell$ in this case.\hfill $\qedbox$

The next lemma presents exactly the family of functions for which $\beta_\ell$
is constant; this technical result will be used in the proof of Theorem \ref{theorem:thm-1.8}.

\begin{lemma}   \label{theorem:solve-beta} Let $\mathcal{O} \subseteq \mathbb{R}$, and denote $\mathcal{O}^p$ as the product of $\mathcal{O}$ with itself $p$ times.
\begin{enumerate}
\item  Let $g_i: \mathcal{O}
\rightarrow
\mathbb{R}$.   Let $g_i\in C^\infty(\mathcal{O})$ for $1\le i\le p$. Suppose
that $\sum_{i=1}^pg_i(u_i)$ is constant on $\mathcal{O}^p$. Then $g_i$ is
constant for $1\le i\le p$.

\item  Suppose $f^{(2)}(0) \neq 0$, and $k \in \mathbb{R}$.  Then the local solutions to the differential equation $\Omega(f) = 
\frac{f^{(3)}(1+f')}{\left[f^{(2)}\right]^2}= k$ are as follows:
\begin{enumerate}
\item  $k = 0 \Rightarrow f$ is quadratic.
\item  $k = 1 \Rightarrow 1+f' = e^{au+b}$ for some $0<a \in \mathbb{R}$, and $b
\in \mathbb{R}$.
\item  $k \neq 0$ and $k \neq 1 \Rightarrow 1+f' = \sqrt[1-k]{(1-k)(au+b)}$ for some
$0<a\in \mathbb{R}$ and $b \in \mathbb{R}$.
\end{enumerate}

\item  Any solution to $\beta_2 = k$  where $k$ is constant is also a solution to
$\beta_\ell = k'$ where $k'$ is constant.
\end{enumerate}
\end{lemma}

\begin{proof}

Assertion 1 is obvious as each summand is a function of different variables.  We apply the previous assertion to the differential equation $\beta_2 = k$ to
note that each of the summands 
$\frac{f_j^{(3)}(1+f_j')}{\left[f_j^{(2)}\right]^2}$ is constant.  We
can solve this explicitly for all functions on which $\beta_\ell$ is defined. 
The hypotheses ensure that the given expression makes sense in a small
neighborhood of $u = 0$.   We consider each case given in the theorem:

 Case I:  $k = 0$.  This is more or less obvious since the denominator  of
$\Omega$ is nonzero, and 
$(1+f')$ is nonzero.  Thus  $f^{(3)} = 0$; this establishes Assertion 2(a).
For the next cases, we compute

\begin{equation}  \label{equation:eqn-5.a}
\begin{array}  {r c l}
\frac{f^{(3)}(1+f')}{\left[f^{(2)}\right]^2} & = & k  \Longleftrightarrow  \\
\frac{f'''}{f''}  & = &  \frac{f''}{1+f'}k    \Longleftrightarrow   \\
\log f''& =& k\log (1+f') + a'  \Longleftrightarrow \\
\frac{f''}{(1+f')^k} & = & e^{a'} = a\,.
\end{array}
\end{equation}
Case II:  $k = 1$.  We integrate Equation (\ref{equation:eqn-5.a}) to get
\begin{equation*}
\begin{array}  {r c l}
\log(1+f') & = &au + b  \Longleftrightarrow  \\
1+f' & = & e^{au+b}\,.
\end{array}
\end{equation*}
Case III:  $k \neq 0$ and $k \neq 1$.  We integrate  (\ref{equation:eqn-5.a}) to get
\begin{equation*}
\begin{array}  {r c l}
\frac{1}{1-k}(1+f')^{1-k}  & = & au + b  \Longleftrightarrow  \\
1+f' & = & \sqrt[1-k]{(1-k)(au+b)}\,.
\end{array}
\end{equation*}
One can simply check that each of the families found in in the previous
assertion are also solutions to $\beta_\ell = $ constant.  Of course, more
initial conditions will need to be given for higher values of $\ell$ to
completely describe all solutions.
\end{proof}

We will need another family of invariants  can be constructed in the same manner as
$\beta_\ell$ using the other nonzero higher covariant
derivatives of the curvature tensor $R$, as listed in Lemma
\ref{theorem:derivatives}.  Here, we may  remove the hypothesis that $f_i''
\neq 0$.

\begin{theorem}  \label{theorem:gammas}  Adopt the notation of Definitions \ref{theorem:def} and \ref{theorem:manifolds}, and let $\BB$ be a normalized basis. Suppose $\ell \geq 2$, and set 
$$
\gamma_\ell = \sum_j 
\nabla^\ell R(U_0, U_j, U_j, U_0; U_j, \ldots, U_j) \cdot \nabla 
R(U_0, U_j, U_j, U_0; U_j)^{\ell-2}\,.
$$
\begin{enumerate}

\item $\gamma_\ell$ is independent of the normalized basis chosen, and is an
$\ell$-model invariant.

\item  $\gamma_2 = \sum_j \left[
\frac{\ep_j}{(f_j'+1)^2}
\left(
4(f_j')^2 +2f_j' + 6f_jf_j'' - \frac{(f_j)^2f_j'''}{f_j'+1}
\right)\right]\,.$

\item  If $\mathcal{M}_F$ is $\ell$-curvature homogeneous, then $\gamma_p$ is
constant for $1\le p \le \ell$.

\item  If $\mathcal{M}_F$ is locally homogeneous, then $\gamma_\ell$ is constant
for all $\ell$.

\end{enumerate}
\end{theorem}

\begin{proof}
Let $\tilde \BB$ be another normalized basis. By Theorem
\ref{theorem:thm-1.5} there exists $a_0 \neq 0$ and a $\sigma \in
\operatorname{Sym}_k$ so that 

\begin{eqnarray*}
&&\nabla^\ell R(\tilde U_0, \tilde U_j, \tilde U_j, \tilde U_0; \tilde U_j, \ldots,
\tilde U_j)\\&  = & \left( \frac{1}{\sqrt{|a_0|}} \right)^{\ell - 2}
\nabla^\ell R(U_0, U_{j'}, U_{j'}, U_0; U_{j'}, \ldots, U_{j'}),\end{eqnarray*}

and

$$
 \nabla R(\tilde U_0, \tilde U_j, \tilde U_j, \tilde U_0; \tilde U_j)  = 
 \sqrt{|a_0|}
\nabla R(U_0, U_{j'}, U_{j'}, U_0; U_{j'})\,.
$$
where $j' = \sigma(j)$.  Combining the above according to the definition of
$\gamma_\ell$ establishes Assertion 1.  Lemma  \ref{theorem:nabla-on-norm} and Theorem \ref{theorem:derivatives} establishes Assertion 2.

Assertions 3 and 4 follow similarly as in  the proof of Assertions
2 and 3 of Theorem \ref{theorem:extra}.
\end{proof}

 We use the invariants described above to study the local homogeneity of the
manifold $\mathcal{M}_F$, and establish Theorem \ref{theorem:thm-1.8}.

\emph{Proof of Theorem \ref{theorem:thm-1.8}.}
If $\mathcal{M}_F$ were 2-curvature homogeneous, then by
Assertion 3 of Theorem \ref{theorem:independent},  $\beta_2$ is
constant.   By Assertion 3 of Theorem \ref{theorem:gammas}, $\gamma_2$ must also be constant.  None of the solutions to $\beta_2 = $ constant listed
in Lemma
\ref{theorem:solve-beta} make $\gamma_2$ constant as well.
\hfill $\qedbox$

In most cases, Theorem \ref{theorem:thm-1.8} tells us these manifolds are not
2-curvature homogeneous, and hence not generally locally homogeneous.  One asks
if any of the $\mathcal{M}_F$ are 1-curvature homogeneous.  We
will study this question in a subsequent paper.

\section*{Acknowledgments}
The author would like to take this opportunity to thank P. Gilkey and E. Puffini for their valuable help while this research was conducted.  Also, the author wishes to thank the referee for several valuable suggestions.  The research of the author is partially supported by a CSUSB faculty research grant.


\begin{thebibliography}{AAA}

\bibitem{BKV}  E. Boeckx, O. Kowalski, and L. Vanhecke, \textbf{Riemannian manifolds of conullity two}, World Scientific, (1996), ISBN: 981-02-2768-X.

\bibitem{Bue}  P. Bueken,  {\it On curvature homogeneous three-dimensional
Lorentzian manifolds}, Journal of Geometry and Physics, {\bf 22}, (1997),
349--362.


\bibitem{BD}  P. Bueken, and M. Djori\'c, {\it Three-dimensional Lorentz metrics and curvature homogeneity of order one,} Ann. Global Anal.
Geom.,
\textbf{18}, (2000), 85--103.

\bibitem{BV}  P. Bueken, and L. Vanhecke, {\it Examples of curvature homogeneous Lorentz metrics,}  Classical Quantum Gravity, \textbf{14},
(1997), L93--L96.

\bibitem{Derd}  A. Derdzinski, {\it Einstein metrics in dimension four}, \textbf{Handbook of differential geometry.  Vol. I.}, Edited by Dillen and Verstaelen, North-Holland, Amsterdam, (2000),  419--707, ISBN: 0-444-82240-2.

\bibitem{DG}  C. Dunn, and  P. Gilkey, {\it Curvature homogeneous pseudo-Riemannian
manifolds which are not locally homogeneous}, \textbf{Complex, Contact and
Symmetric Manifolds}, Birkh$\ddot{\hbox{a}}$user  (2005), 145--152, ISBN:
0-8176-3850-4.


\bibitem{DGS}  C. Dunn, and P. Gilkey,   and S. Nik\v cevi\'c,  {\it Curvature homogeneous
signature $(2,2)$ manifolds}, {\bf Differential Geometry and its
Applications, Proceedings of the 9th International Conference}, (2004), 29--44,
ISBN:  80-86732-63-0.


\bibitem{FKM}  D. Ferus, H. Karcher, and  H. M\"unzner, {\it Cliffordalgebren und neue isoparametrische Hyperfl\"achen}, Math. Z., \textbf{177}, (1981), 479--502.

\bibitem{GBk}  P. Gilkey, \textbf{Geometric Properties of Natural Operators
Defined by the Riemann Curvature Tensor}, World Scienific (2001), ISBN:
981-02-4752-4.

\bibitem{GBk2}  P. Gilkey, \textbf{The Geometry of Curvature Homogeneous
Pseudo--Rie\-mannian Manifolds}, Imperial College Press (2006), to appear.


\bibitem{GIZ}  P. Gilkey, R. Ivanova, and  T. Zhang, {\it Szabo Osserman IP
pseudo-Riemannian manifolds},  Publ. Math. Debrecen, \textbf{62}, (2003),
387-401.

\bibitem{GKV}  E. Garc\'ia Rio, D. Kupeli,  and  R. V\'azquez-Lorenzo, \textbf{Osserman Manifolds in Semi-Riemannian Geometry}, Lecture notes in
Mathematics, Spring-Verlag, (2002), ISBN: 3-540-43144-6.

\bibitem{GSL} P. Gilkey, and S. Nik\v cevi\'c, {\it Affine curvature homogeneous
3-dimensional Lorentz manifolds},  Int. J. Geom. Meth. Mod. Phys., {\bf 2},
(2005), 737--749.

\bibitem{GS}  P. Gilkey,  and  S. Nik\v cevi\'c, {\it Complete curvature homogeneous
pseudo-Riemannian manifolds}, Classical and Quantum
Gravity, {\bf 21}, (2004), no. 15, 3755--3770.

\bibitem{GSs} P. Gilkey,  and S. Nik\v cevi\'c, {\it  Nilpotent spacelike Jordan Osserman pseudo-Riemannian manifolds}, Rend. Circ. Mat.
Palermo, (2), Suppl. No. 72, (2004), 99--105. 

\bibitem{GS04} P. Gilkey, and S. Nik\v cevi\'c,  {\it Complete $k$-curvature
homogeneous pseudo-Riemannian manifolds}, Annals of Global Analysis, {\bf 27}
(2005), 87--100.

\bibitem{GSI} P. Gilkey, and S. Nik\v cevi\'c, \emph{Isometry groups of $k$-curvature homogeneous manifolds}, Rendiconti del Circolo Matematico di Pelermo, Proceedings of 25th Winter School Geometry and Physics, Smi 2005, to appear.  math.DG/0505598.
 \bibitem{G-S}  P. Gilkey,  and S. Nik\v cevi\'c, {\it Generalized plane wave
manifolds}, Kragujevac Journal of Mathematics,  {\bf 28}, (2005), 113--138.



\bibitem{Gr} M. Gromov,
{\bf Partial Differential Relations},
Ergeb. Math. Grenzgeb 3. Folge, Band 9, Springer-Verlag (1986), ISBN:
3-540-12177-3.



\bibitem{KTV} O. Kowalski,   F. Tricerri,  and  L. Vanhecke, \emph{Curvature homogeneous
Riemannian manifolds,} J. Math. Pures Appl. (9), {\bf 71}, (1992), no. 6,
471--501.


\bibitem{Op}  B. Opozda, {\it Affine versions of Singer's theorem on locally
homogeneous spaces}, Ann. Global Anal. Geom.,  \textbf{15}, (1997), 187--199.


\bibitem{PS04} F. Podesta, and A. Spiro, {\it Introduzione ai  Gruppi di
Trasformazioni},  Volume of the Preprint Series of the Mathematics
Department ``V. Volterra'' of the University of Ancona,  Via delle Brecce
Bianche, Ancona, ITALY, (1996).

\bibitem{PTV} F.  Pr\"ufer, F. Tricerri, and L. Vanhecke, {\it Curvature invariants, differential operators and local homogeneity,} Trans. Am. Math. Soc., \textbf{348}, (1996), 4643--4652.

\bibitem{SSV1}  K. Sekegawa, H. Suga,  and  L. Vanhecke, {\it Four-dimensional curvature homogeneous spaces,} Commentat. Math. Univ. Carol.,
\textbf{33}, (1992), 261--268.

\bibitem{SSV2}  K. Sekegawa, H. Suga,  and  L. Vanhecke, {\it Curvature homogeneity for four-dimensional manifolds}, J. Korean Math. Soc.,
\textbf{32}, (1995), 93--101.

\bibitem{S} I. M. Singer, {\it Infinitesimally homogeneous spaces},  Commun.
Pure Appl. Math., {\bf 13}, (1960), 685--697.

\bibitem{T}  H. Takagi, {\it On curvature homogeneity of Riemannian manifolds}, T$\hat{\hbox{o}}$hoku Math. J., \textbf{26}, (1974), 581--585.

\bibitem{TV} F. Tricerri, and L. Vanhecke, {\it Vari\'et\'es riemanniennes dont
le tenseur de courbure est celui d'un espace sym\'etrique riemannien
irr\'eductible},  C. R. Acad. Sci., Paris, S\'er. I, {\bf 302}, (1986),
233--235.

\bibitem{Y} K. Yamato, {\it Algebraic Riemann manifolds,}  Nagoya Math. J., \textbf{115}, (1989), 87--104.


\end{thebibliography}
\end{document}